\documentclass[11pt, reqno]{amsart}
\usepackage{amssymb}
\usepackage{amscd}
\usepackage{verbatim,ifthen}
\usepackage{color}
\usepackage{latexsym}
\usepackage{tikz}
\usepackage{tikz-cd}
\usepackage{mathrsfs}
\usepackage{wrapfig}
\usetikzlibrary{shapes}
\usepackage{color}
\usetikzlibrary{arrows.meta}
\usepackage{bbm}
\usetikzlibrary{matrix}
\usetikzlibrary{calc}
\usetikzlibrary{arrows,intersections}
\usepackage{pgfplots}
\usepackage{multicol}
\usepackage{array}
\newcolumntype{M}[1]{>{\centering\arraybackslash}m{#1}}

\usepackage[colorlinks, linkcolor=black, citecolor=magenta, linktocpage,backref=page]{hyperref}
\renewcommand*{\backref}[1]{}
\renewcommand*{\backrefalt}[4]{[{\tiny%
    \ifcase #1 Not cited.%
          \or Cited on page~#2.%
          \else Cited on pages #2.%
    \fi%
    }]}
\addtocontents{toc}{\setcounter{tocdepth}{1}}
\usepackage{amsmath}

\numberwithin{equation}{section}

\let\oldtocsection=\tocsection
 
\let\oldtocsubsection=\tocsubsection

\renewcommand{\tocsection}[2]{\hspace{0em}\oldtocsection{#1}{#2}}
\renewcommand{\tocsubsection}[2]{\hspace{1em}\oldtocsubsection{#1}{#2}}

\def\XXint#1#2#3{{\setbox0=\hbox{$#1{#2#3}{\int}$ }
\vcenter{\hbox{$#2#3$ }}\kern-.6\wd0}}

\usepackage{eucal}
\usepackage{calc}  
\usepackage{enumitem} 
\usepackage{tensor}
\usepackage{graphicx,wrapfig,lipsum}
\usepackage{etoolbox}
\usepackage{marginnote}
\usepackage{lipsum}
\makeatletter
\patchcmd{\@mn@margintest}{\@tempswafalse}{\@tempswatrue}{}{}
\patchcmd{\@mn@margintest}{\@tempswafalse}{\@tempswatrue}{}{}
\reversemarginpar 
\makeatother
\usepackage{scrextend}

\makeatletter
\DeclareRobustCommand\widecheck[1]{{\mathpalette\@widecheck{#1}}}
\def\@widecheck#1#2{%
    \setbox\z@\hbox{\m@th$#1#2$}%
    \setbox\tw@\hbox{\m@th$#1%
       \widehat{%
          \vrule\@width\z@\@height\ht\z@
          \vrule\@height\z@\@width\wd\z@}$}%
    \dp\tw@-\ht\z@
    \@tempdima\ht\z@ \advance\@tempdima2\ht\tw@ \divide\@tempdima\thr@@
    \setbox\tw@\hbox{%
       \raise\@tempdima\hbox{\scalebox{1}[-1]{\lower\@tempdima\box
\tw@}}}%
    {\ooalign{\box\tw@ \cr \box\z@}}}
\makeatother
\title{The Schwarz Lemma in K\"ahler and Non-K\"ahler Geometry}
\author{Kyle Broder}

\makeatother
\usepackage{xcolor}

\def\XXint#1#2#3{{\setbox0=\hbox{$#1{#2#3}{\int}$ }
\vcenter{\hbox{$#2#3$ }}\kern-.6\wd0}}

\usepackage{soul}

\usepackage{calc}  
\usepackage{enumitem} 
\usepackage{tensor}
\usepackage{graphicx,wrapfig,lipsum}
\usepackage{etoolbox}
\usepackage{marginnote}
\usepackage{lipsum}
\makeatletter
\patchcmd{\@mn@margintest}{\@tempswafalse}{\@tempswatrue}{}{}
\patchcmd{\@mn@margintest}{\@tempswafalse}{\@tempswatrue}{}{}
\reversemarginpar 
\makeatother
\usepackage{array}
\newcolumntype{M}[1]{>{\centering\arraybackslash}m{#1}}

\usepackage{letltxmacro}
\makeatletter
\AtBeginDocument{%
  \@ifdefinable{\myorg@nameref}{%
    \LetLtxMacro\myorg@nameref\nameref
    \DeclareRobustCommand*{\nameref}[1]{%
      \emph{\myorg@nameref{#1}}%
    }%
  }%
}
\makeatother

\theoremstyle{plain}
\theoremstyle{definition}

\thanks{The author was partially supported by an Australian Government Research Training Program (RTP) Scholarship and funding from the Australian Government through the Australian Research Council's Discovery Projects funding scheme (project DP220102530).}

\begin{document}

\maketitle

\begin{abstract}
We introduce a new curvature constraint that provides an analog of the real bisectional curvature considered by Yang--Zheng \cite{YangZhengRBC} for the Aubin--Yau inequality.  A unified perspective of the various forms of the Schwarz lemma is given, leading to novel Schwarz-type inequalities in both the K\"ahler and Hermitian categories.  
\end{abstract}

\section{Introduction}

The holomorphic sectional curvature of a Hermitian metric $\omega$,  $$\text{HSC}_{\omega} : \mathcal{F}_M \times \mathbb{C}^n \backslash \{ 0 \} \longrightarrow \mathbb{R}, \hspace{1cm} \text{HSC}_{\omega}(\vartheta, v) := \frac{1}{| v |_{\omega}^4} \sum_{i,j,k, \ell} R_{i \overline{j} k \overline{\ell}} v_i \overline{v}_j v_k \overline{v}_{\ell},$$ is one of the most curious objects in complex geometry. This supernatural object determines the curvature tensor entirely if the metric is K\"ahler, but is not strong enough to dominate\footnote{We say that a curvature $P$ dominates a curvature $Q$ if a bound on $P$ forces a similar bound on $Q$.} the Ricci curvature. The relationship between the holomorphic sectional curvature and the Ricci curvature(s) has been the subject of great interest for decades. In the hierarchy of curvatures in complex geometry, they occupy a shared class: They are both dominated by the holomorphic bisectional curvature, and both dominate the scalar curvature\footnote{The holomorphic sectional curvature of a \textit{Hermitian metric} only dominates the sum $\text{Scal}_{\omega} + \widetilde{\text{Scal}}_{\omega}$ of scalar curvature $\text{Scal}_{\omega} : = \text{tr}_{\omega} \text{Ric}_{\omega}^{(1)}$ and the altered scalar curvature $\widetilde{\text{Scal}}_{\omega} : = \text{tr}_{\omega} \text{Ric}_{\omega}^{(3)}$.  These scalar curvatures coincide if the metric is K\"ahler (or more generally, if the metric is balanced).}. But do not dominate each other: Hitchin's examples \cite{Hitchin} of Hodge metrics on Hirzebruch surfaces have $\text{HSC} >0$, but do not admit K\"ahler metrics with $\text{Ric} >0$.  On the other hand,  there are Fermat hypersurfaces of degree $d \geq n+2$ in $\mathbb{P}^n$ which have K\"ahler metrics with $\text{Ric} <0$ but no Hermitian metrics with $\text{HSC}<0$ (see, e.g., \cite{WuYau1}).

It goes back to Grauert--Reckziegel \cite{GrauertReckziegel} that a compact complex manifold $M$ supporting a Hermitian metric of negative holomorphic sectional curvature is Kobayashi hyperbolic, i.e., every entire curve $\mathbb{C} \to M$ is constant (see \cite{GreeneWu} for more refined statements).  The converse is false; an example was constructed by Demailly in \cite[Theorem 8.2]{DemaillyExample} (c.f., \cite{DiverioSurvey}). The negativity of the holomorphic sectional curvature describes hyperbolicity in terms of lengths, while the negativity of the Ricci curvature expresses hyperbolicity in terms of volumes. The heuristic that lengths control volumes leads one to recover the famous conjecture made by Yau some decades ago:\begin{eqnarray}\label{YauConjecture}
\text{HSC}_{\omega} < 0 & \implies & K_M > 0.
\end{eqnarray}

Of course, if $(M,\omega)$ is K\"ahler, Yau's solution of the Calabi conjecture \cite{Yau} implies that $K_M>0$ is equivalent to the existence of an (a priori different) K\"ahler metric with negative Ricci curvature.  Yau's conjecture \eqref{YauConjecture} was verified for K\"ahler surfaces by Campana \cite{Campana} and Wong \cite{Wong},  making use of the Enriques--Kodaira classification.  The major breakthrough which incited the tremendous progress that followed was achieved by Heier--Lu--Wong \cite{HeierLuWong}, proving \eqref{YauConjecture} for projective threefolds.  This was extended to all dimensions by (Wong--)Wu--Yau in \cite{WongWuYau, WuYau1,WuYau2} and subsequently extended to the K\"ahler setting by Tosatti--Yang \cite{TosattiYang}. Diverio--Trapani \cite{DiverioTrapani} later relaxed the negativity of the holomorphic sectional curvature to quasi-negativity (i.e., non-positive everywhere and negative at one point).   It is natural to ask whether the same result holds if the metric is merely Hermitian; that is, suppose a compact complex manifold admits a Hermitian metric of negative sectional curvature, is the canonical bundle $K_M$ ample? 

Once the smoke has cleared, this question levels to a problem concerning the Schwarz lemma. That is, understanding the conditions which force control on the growth of holomorphic maps. More precisely, a lemma is of Schwarz-type if, for any holomorphic map $f :(M,\omega) \longrightarrow (N, \eta)$ between Hermitian manifolds, there are assumptions depending only on the curvature of the metrics $\omega$ and $\eta$, and geometric quantities such as rank and dimension, such that the energy density $| \partial f |^2$ of $f$ is bounded by some function of these curvature assumptions. 

There are two leading Schwarz lemmas: The Chern--Lu inequality estimates $| \partial f |^2$ if one has control of the Ricci curvature of $\omega$ from below, and the holomorphic sectional curvature of $\eta$ from above. The Aubin--Yau inequality estimates $| \partial f |^2$ if one has control of the Ricci curvature of $\eta$ from above, and the holomorphic sectional curvature of $\omega$ from below.

In \cite{RubinsteinKE}, Rubinstein detailed a framework that unified these two front-runners. The underlying philosophy is that these theorems are best understood via holomorphic maps, and not as abstract tensor calculations (c.f., \cite{BroderSurvey}). This paradigm extends to the study of holomorphic maps between Hermitian manifolds and can be made further transparent by using the ideas in \cite{YangZhengRBC}.

Indeed, it goes back to Ahlfors \cite{Ahlfors} that an estimate on $|\partial f |^2$ is obtained from applying the maximum principle. Hence, we require a Laplace operator. Two such operators are given to us: the complex Laplacians $\Delta_{\omega} : = \text{tr}_{\omega} \sqrt{-1} \partial \overline{\partial}$ and $\Delta_{\eta} := \text{tr}_{\eta} \sqrt{-1} \partial \overline{\partial}$. The Bochner formula \cite{KobayashiWu} details how to compute the complex Hessian of (the norm of) a section $\sigma \in H^0(\mathcal{E})$ of a holomorphic vector bundle $\mathcal{E}$: \begin{eqnarray*}
\partial \overline{\partial} | \sigma |^2 &=& \langle \nabla \sigma, \nabla \sigma \rangle -    \langle R^{\mathcal{E}} \sigma, \sigma \rangle,
\end{eqnarray*}

where $R^{\mathcal{E}}$ is the curvature of the Chern connection on $\mathcal{E}$. When $\sigma = \partial f$, the holomorphic vector bundle $\mathcal{E}$ is the twisted cotangent bundle $\mathcal{T}_M^{\ast} \otimes f^{\ast} \mathcal{T}_N$. Since the curvature of the tensor product of bundles splits additively, we get opposing contributions to the curvature from the source and target metrics. The trace of the Bochner formula with respect to the source metric yields the Chern--Lu inequality; the trace of the Bochner formula with respect to the target metric yields the Aubin--Yau inequality.  Note that this interpretation of the Aubin--Yau second-order estimate was Yau's motivation in 1976.

The trace of the curvature with respect to the same metric (i.e., the trace of $R^{\mathcal{T}_M}$ with respect to $\omega$, or the trace of $R^{\mathcal{T}_N}$ with respect to $\eta$) is easy to understand and control: it is the second Chern--Ricci curvature. The trace of the curvature with respect to the opposing metric, however, is much more mysterious. For K\"ahler metrics, Yau \cite{YauSchwarzLemma} controlled this curvature term by assuming control of the bisectional curvature. Royden \cite{Royden} refined this, in the K\"ahler setting, to bounds on the holomorphic sectional curvature. 

If the target metric does not happen to be a God-fearing K\"ahler metric, the target curvature term arising in the Schwarz lemma does not appear to be controlled by the holomorphic sectional curvature. This motivated the real bisectional curvature that was introduced by Yang--Zheng \cite{YangZhengRBC}: $$\text{RBC}_{\omega} : \mathcal{F}_M \times \mathbb{R}_+^n  \longrightarrow \mathbb{R}, \hspace*{1cm} \text{RBC}_{\omega}(\vartheta, \lambda) := \frac{1}{| \lambda |^2} \sum_{\alpha, \gamma} R_{\alpha \overline{\alpha} \gamma \overline{\gamma}} \lambda_{\alpha} \lambda_{\gamma}.$$ Here, $\mathcal{F}_M$ denotes the unitary frame bundle, $R_{\alpha \overline{\beta} \gamma \overline{\delta}}$ denote the components of the Chern curvature tensor with respect to the frame $\vartheta$, and $\lambda = (\lambda_1, ..., \lambda_n) \in \mathbb{R}_+^n = \{ (v_1, ..., v_n) \in \mathbb{R}^n \backslash \{ 0 \} : v_k \geq 0, \forall k \}$.  With this refined curvature constraint,  we have the following Chern--Lu inequality due to Yang--Zheng \cite{YangZhengRBC}:

\subsection*{Theorem 1.1}\label{HermitianCL}
(Hermitian Chern--Lu). Let $f : (M^n, \omega) \longrightarrow (N, \eta)$ be a holomorphic map between Hermitian manifolds. Assume $\text{Ric}_{\omega}^{(2)}  \geq - C_1 \omega + C_2 f^{\ast} \eta$ for continuous functions $C_1$ and $C_2 \geq 0$. Assume that $\text{RBC}_{\eta} \leq - \kappa \leq  0$ for some continuous function $\kappa$. Then \begin{eqnarray*}
\Delta_{\omega} \log |\partial f |^2 & \geq & - C_1 + \frac{1}{r} \left(   f^{\ast} \kappa  + C_2 \right) |\partial f |^2,
\end{eqnarray*}

In particular, if $\kappa, C_1, C_2$ are constants such that $\kappa + C_2 > 0$, and $M$ is compact, then \begin{eqnarray*}
| \partial f |^2 & \leq &  \frac{C_1 r}{ \kappa + C_2}.
\end{eqnarray*}

The above Schwarz lemma was used by Yang--Zheng \cite{YangZhengRBC} to extend the Wu--Yau theorem, proving that if a compact K\"ahler manifold admits a Hermitian metric with negative real bisectional curvature, then the canonical bundle is ample.

The main object of the present paper is to embed the real bisectional curvature into an algebraic framework -- interpreting the real bisectional curvature as a quadratic form-valued map on the unitary frame bundle.  In more detail,  for any fixed frame, we let $\mathcal{R} \in \mathbb{R}^{n \times n}$ be the matrix with entries $\mathcal{R}_{\alpha \gamma} : = R_{\alpha \overline{\alpha} \gamma \overline{\gamma}}$.  For $v \in \mathbb{R}_+^n$,  the real bisectional curvature is then realized as the Rayleigh quotient $$\text{RBC}_{\omega}(v) \ = \ \frac{v^t \mathcal{R} v}{v^t v}.$$  From this framework,  we not only gain deeper insight into the real bisectional curvature but discover the appropriate analog of the real bisectional curvature for the Aubin--Yau inequality in the Hermitian category: Let $\Gamma_{\circ} : = \{ (v_1, ..., v_n) \in \mathbb{R}_+^n : v_1 \geq v_2 \geq \cdots \geq v_n > 0 \}$.  For $v = (v_1, ..., v_n) \in \Gamma_{\circ}$,  we denote by $u_v : = v_{\circ}^{-1} : = (v_1^{-1}, ..., v_n^{-1})$ the vector obtained from inverting $v$ with respect to the Hadamard product.  The curvature term appearing in the Aubin--Yau inequality is precisely the following:

\subsection*{Definition 1.2}
Let $(M,\omega)$ be a Hermitian manifold.  The \textit{Schwarz bisectional curvature} is the function $$\text{SBC}_{\omega} : \mathcal{F}_M \times \Gamma_{\circ} \longrightarrow \mathbb{R}, \hspace*{1cm} \text{SBC}_{\omega}(\vartheta, v) \ : = \ u_v^t \mathcal{R} v,$$ where $\vartheta$ is a section of the unitary frame bundle $\mathcal{F}_M$. We understand bounds on the Schwarz bisectional curvature as bounds on $u_v^t \mathcal{R} v$ with respect to \textit{all frames}.

\hfill

The corresponding Aubin--Yau inequality appears to be new (even in the K\"ahler case):

\subsection*{Theorem 1.3}\label{HermitianAubinYau}
(Hermitian Aubin--Yau). Let $f : (M, \omega) \longrightarrow (N, \eta)$ be a holomorphic map between Hermitian manifolds, which is biholomorphic onto its image. Assume $\text{SBC}_{\omega} \geq - \kappa$ and $\text{Ric}_{\eta}^{(2)} \leq - C_1 \eta + C_2 (f^{-1})^{\ast} \omega$ for some continuous functions $\kappa, C_1, C_2$, with $\kappa \geq 0$. Then \begin{eqnarray*}
\Delta_{\eta} | \partial f |^2 & \geq & C_1 | \partial f |^2 - n( C_2 + \kappa).
\end{eqnarray*}

In particular, if $\kappa, C_1, C_2$ are constants such that $C_1 >0$, and $M$ is compact,  \begin{eqnarray*}
| \partial f |^2 & \leq & \frac{n(C_2 + \kappa)}{C_1}.
\end{eqnarray*}

\hfill

One can combine these Schwarz lemmas to obtain the following $8$--dimensional family of Schwarz lemmas:

\subsection*{Theorem 1.4}\label{ParametrizedSL}
Let $f : (M^n, \omega) \longrightarrow (N, \eta)$ be a holomorphic map of rank $r$ between Hermitian manifolds with $\text{SBC}_{\omega} \geq - \kappa_1$ and $\text{RBC}_{\eta} \leq - \kappa_2$, for some constants $\kappa_1, \kappa_2 \geq 0$. Assume there is a Hermitian metric $\mu$ on $M$ such that, for constants $C_1, C_2, C_3, C_4 \in \mathbb{R}$, with $C_2 \geq 0$,  $C_3>0$,  and $C_2 + \kappa_2 >0$, we have \begin{eqnarray*}
- C_1 \mu + C_2 f^{\ast} \eta \ \leq \ \text{Ric}^{(2)}_{\mu} \ \leq \ - C_3 \mu + C_4 \omega.
\end{eqnarray*}

Then, if $M$ is compact, \begin{eqnarray*}
| \partial f |^2 & \leq & \frac{C_1 n r(\kappa_1 + C_4)}{C_3(\kappa_2 + C_2)}.
\end{eqnarray*}

\hfill

One particular corollary of the above theorem is the following Hermitian analog of the Chen--Cheng--Lu Schwarz lemma  \cite{ChenChengLu}: 

\subsection*{Corollary 1.5}
Let $f : (M^n, \omega) \longrightarrow (N, \eta)$ be a holomorphic map of rank $r$ between Hermitian manifolds with $\text{SBC}_{\omega} \geq - \kappa_1$ and $\text{RBC}_{\eta} \leq - \kappa_2$, for some constants $\kappa_1, \kappa_2 \geq 0$. Assume there is a Hermitian metric $\mu$ on $M$ such that $$-C_1 \mu + C_2 f^{\ast} \eta \ \leq \  \text{Ric}_{\mu}^{(2)} \ \leq \ - C_3 \mu,$$ where $C_1 = \frac{\kappa_2 + C_2}{\kappa_2 nr} C_3$, and $C_2 \geq \kappa_2(nr-1)$. Then, if $M$ is compact,  \begin{eqnarray*}
| \partial f |^2 & \leq & \frac{\kappa_1}{\kappa_2}.
\end{eqnarray*}

The theorem also yields a Schwarz lemma expressed exclusively in terms of (second) Ricci curvatures: 

\subsection*{Corollary 1.6}
Let $f : (M^n, \omega) \longrightarrow (N, \eta)$ be a holomorphic map of rank $r$ between Hermitian manifolds with $\text{SBC}_{\omega} \geq - \kappa_1$ and $\text{RBC}_{\eta} \leq - \kappa_2$, for some constants $\kappa_1, \kappa_2 \geq 0$. Assume there is a Hermitian metric $\mu$ on $M$ such that $$-C_1 \mu + C_2 f^{\ast} \eta \ \leq \  \text{Ric}_{\mu}^{(2)} \ \leq \ - C_3 \mu + C_4 \omega,$$ where $C_2 \geq 0$, $C_3 >0$ and $nr(\kappa_1 + C_4) \leq \kappa_2 + C_2$. Then, if $M$ is compact, $$| \partial f |^2 \ \leq \ \frac{C_1}{C_3}.$$

\subsection*{Corollary 1.7}
Let $M^n$ be a compact Hermitian manifold supporting a Hermitian metric $\omega$ with $\text{SBC}_{\omega} \geq - \kappa_1$,and another Hermitian metric $\mu$ such that $$ - C_1 \mu  \ \leq \ \text{Ric}^{(2)}_{\mu} \   \leq \ - C_3 \mu +C_4 \omega ,$$ where $C_1, C_3 >0$ and $C_4 < 0$. If $0 \leq \kappa_1 \leq - C_4$, then there is no holomorphic map into a Hermitian manifold supporting a Hermitian metric $\eta$ with $\text{RBC}_{\eta} <0$.

\subsection*{Remark 1.8}
After releasing an earlier version of this manuscript to the arXiv,  there have been a number of developments that were directly incited by the results of the present manuscript.  For instance,  the paradigm of viewing these curvatures as quadratic form-valued maps on the unitary frame bundle lead to the systematic investigation of the quadratic orthogonal bisectional curvature which appears in \cite{BroderQOBC}.  Surprisingly,  this study, which was motivated entirely within the context of complex differential geometry,  answered questions in graph theory and combinatorics \cite{BroderGraph} and distance geometry \cite{BroderLinearAlgebra}. The considerable amount of attention received from the initial upload of this manuscript was accompanied by encouragement to expand upon the ideas which led to the results presented here.  The details of this narrative were expounded upon in the expository paper \cite{BroderSurvey}.

\subsection*{Acknowledgments}
This work owes much to \cite{RubinsteinKE} and \cite{YangZhengRBC}, from which I have learned a great deal. I'm also grateful to my supervisors Ben Andrews and Gang Tian for their unwavering encouragement and valuable discussions. I would like to thank Fangyang Zheng,  Jeffrey Streets, Kai Tang, Shin-ichi Matsumara, Haojie Chen, and Xiaolan Nie for their interest in this work, and comments on the previous version of this manuscript.  Let me also thank the referee for their detailed and valuable comments.

\section{The Schwarz bisectional curvature}
Let $(M, g,J)$ be a Hermitian manifold. The Chern connection $\nabla$ is the unique connection on $\mathcal{T}_M$ such that $\nabla g = \nabla J =0$ and whose torsion has vanishing $(1,1)$ part. The curvature of the Chern connection in a local coordinate frame $\frac{\partial}{\partial z_i}$ reads \begin{eqnarray*}
R_{i \overline{j} k \overline{\ell}} &=& -\frac{\partial^2 g_{k \overline{\ell}}}{\partial z_i \partial \overline{z}_j} + g^{p \overline{q}} \frac{\partial g_{k \overline{q}}}{\partial z_i} \frac{\partial g_{p \overline{\ell}}}{\partial \overline{z}_j}.
\end{eqnarray*}

To motivate the definition of the Schwarz bisectional curvature, recall that for a holomorphic map $f : (M, \omega) \longrightarrow (N,\eta)$ between Hermitian manifolds, the differential $ \partial f $ is a section of the twisted cotangent bundle $\mathcal{T}_M^{\ast} \otimes f^{\ast} \mathcal{T}_N$. The curvature of the connection induced by the Chern connections on each factor splits according to \begin{eqnarray}\label{SplitMotivation}
R^{\mathcal{T}_M^{\ast} \otimes f^{\ast} \mathcal{T}_N} &=& - R^{\mathcal{T}_M} \otimes \text{id} + \text{id} \otimes f^{\ast} R^{\mathcal{T}_N}.
\end{eqnarray}

The Bochner formula informs us on how to compute the complex Hessian of the energy density: \begin{eqnarray}\label{BochnerMotivation}
\partial \overline{\partial} | \partial f |^2 & = & \langle \nabla^{1,0} \partial f, \nabla^{1,0} \partial f \rangle -  \langle R^{\mathcal{T}^{\ast}_M \otimes f^{\ast} \mathcal{T}_N} \partial f, \partial f \rangle.
\end{eqnarray}

The trace of \eqref{BochnerMotivation} with respect to the source metric $\omega$ furnishes the Chern--Lu inequality. If $f$ is biholomorphic onto its image, the trace of \eqref{BochnerMotivation} with respect to the target metric $\eta$ furnishes the Aubin--Yau inequality. Because of \eqref{SplitMotivation}, the Chern--Lu and Aubin--Yau inequalities require an understanding of the curvatures which arise from the traces of $R^{\mathcal{T}_M} \otimes \text{id}$ and $\text{id} \otimes f^{\ast} R^{\mathcal{T}_N}$ with respect to the metrics on the source or target. 

Taking the trace of \eqref{BochnerMotivation} with respect to $\omega$, we see that in coordinates, we have \begin{eqnarray*}
g^{i \overline{j}} \partial_i \partial_{\overline{j}} (g^{k \overline{\ell}} h_{\gamma \overline{\delta}} f_k^{\gamma} \overline{f_{\ell}^{\delta}}) &=& g^{i \overline{j}} g^{k \overline{\ell}} h_{\gamma \overline{\delta}} f_{ik}^{\gamma} \overline{f_{j\ell}^{\delta}} + g^{i \overline{j}} R_{i \overline{j} k \overline{\ell} }^g g^{k \overline{q}} g^{p \overline{\ell}} h_{\alpha \overline{\beta}} f_p^{\alpha} \overline{f_q^{\beta}} \nonumber \\
&& \hspace*{4cm} - g^{i \overline{j}} g^{p \overline{q}} R_{ \alpha \overline{\beta}\gamma \overline{\delta}}^h f_i^{\alpha} \overline{f_{j}^{\beta}} f_p^{\gamma} \overline{f_{q}^{\delta}},
\end{eqnarray*}

where $g,h$ are the respective metrics underlying $\omega, \eta$. The first term on the right-hand side is, in the language of harmonic maps (c.f., \cite{HarmonicMaps}), the second fundamental form of $f$, namely $| \nabla \partial f |^2$. As indicated in \nameref{SecondRicci}, the second term $g^{i \overline{j}} R_{ i \overline{j}k \overline{\ell}} = \text{Ric}_{k \overline{\ell}}^{(2)}$ is the second Ricci curvature. We feel that we understand this. We want to understand the last term: \begin{eqnarray}\label{MOTIVATIONCURVATURE}
g^{i \overline{j}} g^{p \overline{q}} R_{ \alpha \overline{\beta}\gamma \overline{\delta}}^h f_i^{\alpha} \overline{f_{j}^{\beta}} f_p^{\gamma} \overline{f_{q}^{\delta}}.
\end{eqnarray}

To this end, choose coordinates $(z_1, ..., z_n)$ centered at a point $p \in M$ and $(w_1, ..., w_n)$ at $f(p) \in N$ such that $g_{i \overline{j}} =\delta_{ij}$ and $h_{\alpha \overline{\beta}} = \delta_{\alpha\beta}$ at $p$ and $f(p)$, respectively. If $f = (f^1, ..., f^n)$, then with $f_i^{\alpha} = \frac{\partial f^{\alpha} }{\partial z_i}$, the coordinates can be chosen such that $f_i^{\alpha} = \lambda_i \delta_i^{\alpha}$, where $\lambda_1 \geq \lambda_2 \geq \cdots \geq \lambda_r \geq \lambda_{r+1} = \cdots = 0$, and $r$ is the rank of $\partial f = (f_i^{\alpha})$. Hence, \eqref{MOTIVATIONCURVATURE} reads \begin{eqnarray}\label{RBCFIND}
g^{i \overline{j}} g^{p \overline{q}} R_{ \alpha \overline{\beta}\gamma \overline{\delta}}^h f_i^{\alpha} \overline{f_{j}^{\beta}} f_p^{\gamma} \overline{f_{q}^{\delta}} &=& \sum_{\alpha, \gamma} R^h_{\alpha \overline{\alpha} \gamma \overline{\gamma}} \lambda_{\alpha}^2 \lambda_{\gamma}^2.
\end{eqnarray}

This motivated Yang--Zheng \cite{YangZhengRBC} to consider the following:

\subsection*{Definition 2.1}\label{RBCDefinition}
Let $(M,\omega)$ be a Hermitian manifold. The \textit{real bisectional curvature} $\text{RBC}_{\omega}$ is the function $$\text{RBC}_{\omega} : \mathcal{F}_M \times \mathbb{R}^n \backslash \{ 0 \} \longrightarrow \mathbb{R}, \hspace*{1cm} \text{RBC}_{\omega}(\vartheta, \lambda) := \frac{1}{| \lambda |^2} \sum_{\alpha, \gamma} R_{\alpha \overline{\alpha} \gamma \overline{\gamma}} \lambda_{\alpha} \lambda_{\gamma}.$$ Here, $\mathcal{F}_M$ denotes the unitary frame bundle, $R_{\alpha \overline{\beta} \gamma \overline{\delta}}$ denote the components of the Chern curvature tensor with respect to the frame $\vartheta$, and $\lambda = (\lambda_1, ..., \lambda_n) \in \mathbb{R}_+^n = \{ (v_1, ..., v_n) \in \mathbb{R}^n \backslash \{ 0 \} : v_k \geq 0, \forall k \}$. We say that $\text{RBC}_{\omega} \leq \kappa$ if $\max_{(e,\lambda) \in \mathcal{F}_M \times \mathbb{R}^n \backslash \{ 0 \} } \text{RBC}_{\omega}(e,\lambda) \leq \kappa$. Similar definitions apply for $\text{RBC}_{\omega} \geq \kappa$ and $\text{RBC}_{\omega} = \kappa$. \\

For those readers not so familiar with the real bisectional curvature, let us mention that if the metric is K\"ahler\footnote{Or more generally, K\"ahler-like, i.e., the curvature of the Chern connection has the symmetries of the K\"ahler curvature tensor (see \cite{YangZhengCurvature}).} then the real bisectional curvature is comparable to the holomorphic sectional curvature. For a general Hermitian metric, however, the real bisectional curvature is strictly stronger than the holomorphic sectional curvature. A local example is given in \cite[Example 2.2]{YangZhengRBC}. The real bisectional curvature is not so strong, however, as to control the Ricci curvatures (see \cite[Example 2.3]{YangZhengRBC}). A further discussion of the relationship between the real bisectional curvature and other curvatures is given toward the end of this section.\\

Let us observe, however, that the real bisectional curvature does not give suitable control in the Aubin--Yau second-order estimate. Indeed, assume $f:(M, \omega) \longrightarrow (N, \eta)$ is biholomorphic onto its image. Taking the trace of \eqref{BochnerMotivation}, with respect to $\eta$, we see that, in coordinates,  \begin{eqnarray*}
h^{\gamma \overline{\delta}} \partial_{\gamma} \partial_{\overline{\delta}} \left( g^{i \overline{j}} h_{\alpha \overline{\beta}} f_i^{\alpha} \overline{f_j^{\beta}} \right) &=& h^{\gamma \overline{\delta}} g^{i \overline{j}} h_{\alpha \overline{\beta}} f_{ik}^{\alpha} \overline{f_{j\ell}^{\beta}} (f^{-1})_{\gamma}^k \overline{(f^{-1})_{\delta}^{\ell}} - h^{\gamma \overline{\delta}} g^{i \overline{j}} R^h_{\gamma \overline{\delta} \alpha \overline{\beta}} f_i^{\alpha} \overline{f_j^{\beta}} \\
&& \hspace*{3cm} + h^{\gamma \overline{\delta}} g^{i \overline{q}} g^{p \overline{j}} R^g_{k \overline{\ell} p \overline{q}} h_{\alpha \overline{\beta}} f_i^{\alpha} \overline{f_j^{\beta}} (f^{-1})_{\gamma}^k \overline{(f^{-1})_{\delta}^{\ell}}.
\end{eqnarray*}

As before, the first term on the right-hand side is the pointwise norm squared of the second fundamental form of $f$, and the second term is (minus the) second Ricci curvature of $h$. This time, we want to understand \begin{eqnarray*}
h^{\gamma \overline{\delta}} g^{i \overline{q}} g^{p \overline{j}} R^g_{k \overline{\ell} p \overline{q}} h_{\alpha \overline{\beta}} f_i^{\alpha} \overline{f_j^{\beta}} (f^{-1})_{\gamma}^k \overline{(f^{-1})_{\delta}^{\ell}}.
\end{eqnarray*}

Again, choose coordinates at $p$ and $f(p)$ such that $g_{i \overline{j}} =\delta_{ij}$, $h_{\alpha \overline{\beta}} = \delta_{\alpha\beta}$, and $f_i^{\alpha} = \lambda_i \delta_i^{\alpha}$. Then \begin{eqnarray}\label{AYFIND}
h^{\gamma \overline{\delta}} g^{i \overline{q}} g^{p \overline{j}} R^g_{k \overline{\ell} p \overline{q}} h_{\alpha \overline{\beta}} f_i^{\alpha} \overline{f_j^{\beta}} (f^{-1})_{\gamma}^k \overline{(f^{-1})_{\delta}^{\ell}} &=& \sum_{i,k} R^g_{i \overline{i} k \overline{k}} \lambda_i^2 \lambda_k^{-2}.
\end{eqnarray}

This is not controlled by the real bisectional curvature (what is the vector here?).  \\

To lift the veil on both \eqref{RBCFIND} and \eqref{AYFIND},  introduce the matrix $\mathcal{R} \in \mathbb{R}^{n \times n}$ with entries $\mathcal{R}_{\alpha \gamma} : = R_{\alpha \overline{\alpha} \gamma \overline{\gamma}}$. If $v = (v_1, ..., v_n) \in \mathbb{R}_+^n$,   the real bisectional curvature can be interpreted as the Rayleigh quotient \begin{eqnarray*}
\frac{v^t \mathcal{R} v}{v^t v}.
\end{eqnarray*}

Since the vector $v$ is contained in $\mathbb{R}_+^n$, and not all of $\mathbb{R}^n \backslash \{ 0 \}$,  the (say) negativity of the real bisectional curvature does not precisely correspond to the negative-definiteness of $\mathcal{R}$ (in all frames),  but to the $\mathbb{R}_+^n$--conegativity of $\mathcal{R}$ (in all frames). There is a rich literature concerning this generalization of positive and negative definiteness in convex optimization and control theory. 

This Rayleigh quotient paradigm furnishes insight into the curvature term appearing in \eqref{AYFIND}. Let $\Gamma_{\circ} :=\{(x_1, ..., x_n) \in \mathbb{R}_+^n : x_1 \geq x_2 \geq \cdots \geq x_n > 0 \}$ denote the cone of ordered positive $n$-tuples. For a vector $v \in \Gamma_{\circ}$, we denote by $u_v := v_{\circ}^{-1}$ the vector which inverts $v$ with respect to the Hadamard product. That is, if $v = (v_1, ..., v_n) \in \Gamma_{\circ}$, then $u_v = (v_1^{-1}, ..., v_n^{-1})$.

Then a bound on \eqref{AYFIND} translates to a bound on the \textit{generalized Rayleigh quotient} $$  u_v^t \mathcal{R} v, \hspace{1cm} v \in \Gamma_{\circ}.$$

\subsection*{Remark 2.2}
This yields an interesting comparison between the Chern--Lu and Aubin--Yau inequalities. Indeed, the Chern--Lu inequality requires an upper bound on the real bisectional curvature. The real bisectional curvature is a Rayleigh quotient, which is well-known to give a variational characterization of the eigenvalues.\footnote{Of course,  we are restricting the input arguments of the Rayleigh quotient to the positive orthant,  and hence,  will fail to achieve the eigenvalues, in general. This sentence should be understood heuristically.} The Aubin--Yau inequality requires a lower bound on the Schwarz bisectional curvature. The Schwarz bisectional curvature is a generalized Rayleigh quotient, which is known to give a variational characterization of the singular values. Therefore, at least philosophically, it appears that the Chern--Lu inequality is to the Aubin--Yau inequality what the eigenvalue decomposition is to the singular value decomposition.

\subsection*{Relationship to the Holomorphic Sectional Curvature}

Since we understand the holomorphic sectional curvature a little more than the real bisectional curvature, it is natural to ask when the holomorphic sectional curvature dominates the real bisectional curvature. Of course, in the K\"ahler (or more generally, the K\"ahler-like \cite{YangZhengCurvature}) category, this is true, but no results outside of this obvious case have been considered.

For a non-zero vector $v \in \mathcal{T}_M^{1,0}$, the holomorphic sectional curvature of a Hermitian metric $\omega$ in the direction of $v = (v_1, ..., v_n) \in \mathbb{C}^n \backslash \{ 0 \}$ is defined: \begin{eqnarray*}
\text{HSC}_{\omega}(v) & : = & \frac{1}{| v |_{\omega}^4} \sum_{i,j,k,\ell=1}^n R_{i \overline{j} k \overline{\ell}} v_i \overline{v}_j v_k \overline{v}_{\ell},
\end{eqnarray*}

where $R_{i \overline{j} k \overline{\ell}}$ denotes the components of the Chern curvature tensor in a fixed frame.

It will be convenient to give another description of the holomorphic sectional curvature (in \cite{BroderTangAltered,BroderStanfield}, this is referred to as the \textit{altered holomorphic sectional curvature}). Let $\omega_{\text{FS}}$ denote the Fubini--Study metric on $\mathbb{P}^{n-1}$ of unit volume, and let $[w_1 : \cdots : w_n]$ be the standard unitary homogeneous coordinates. Then \begin{eqnarray*}
\int_{\mathbb{P}^{n-1}} \frac{w_i \overline{w_j} w_k \overline{w_{\ell}}}{| w |^4} \omega_{\text{FS}}^{n-1} &=& \frac{\delta_{ij} \delta_{k\ell} + \delta_{i\ell} \delta_{kj}}{n(n+1)}.
\end{eqnarray*}
Hence, if we fix a point $p \in M$ and any non-negative constants $b_1, ..., b_n$, not all zero, then \begin{eqnarray}\label{HSCCOMPAREFORM}
\sum_{i,j,k,\ell=1}^n \int_{\mathbb{P}^{n-1}} R_{ i\overline{j} k \overline{\ell}} \frac{b_i w_i b_j \overline{w_j} b_k w_k b_{\ell} \overline{w_{\ell}}}{| w |^4} \omega_{\text{FS}}^{n-1} &=& \frac{2}{n(n+1)} \sum_{i,k=1}^n (R_{i \overline{i} k \overline{k}} + R_{i \overline{k} k \overline{i}}) b_i^2 b_k^2.
\end{eqnarray}

In particular, the sign of the holomorphic sectional curvature is controlled by the curvature term $\sum_{i,k=1}^n (R_{i \overline{i} k \overline{k}} + R_{i \overline{k} k \overline{i}}) b_i^2 b_k^2$. That is, $\text{HSC}_{\omega} < 0$ if and only if, for all vectors $v = (v_1, ..., v_n) \in \mathbb{R}_+^n$, we have $$\sum_{i,k} (R_{i \overline{i} k \overline{k}} + R_{i \overline{k} k \overline{i}}) v_i v_k < 0,$$ in any unitary frame. It is clear that $\text{HSC}_{\omega} >0$, $\text{HSC}_{\omega}=0$, etc., can be formulated equivalently, with the obvious changes.

Again, the algebraic framework allows us to make some not-entirely trivial observations:

\subsection*{Theorem 2.3}
Let $(M, \omega)$ be a Hermitian manifold. Suppose the curvature tensor of the Chern connection satisfies \begin{eqnarray}\label{CurvatureSymmetry}
R(X, \overline{X}, Y, \overline{Y}) &=& - R(Y, \overline{Y}, X, \overline{X}),
\end{eqnarray} for all distinct pairs of $(1,0)$-tangent vectors $X,Y \in \mathcal{T}_M^{1,0}$. Then the holomorphic sectional curvature coincides with the real bisectional curvature (up to a positive multiple).\begin{proof}
Suppose $R$ has the above symmetry condition. Expanding $R(X+\lambda Y, \overline{X+\lambda Y}, Z+\mu W, \overline{Z+\mu W}) + R(Z+\mu W, \overline{Z+\mu W}, X+\lambda Y, \overline{X+\lambda Y})=0$ for $\lambda, \mu \in \mathbb{C}$, and looking at the $\mu \overline{\lambda}$ coefficient shows that $$R(X, \overline{Y}, Z, \overline{W}) = - R(Z, \overline{W}, X, \overline{Y}),$$ for all $X,Y,Z,W \in \mathcal{T}_M^{1,0}$ (not all equal). In any fixed frame, the holomorphic sectional curvature is (up to a positive multiple) the Rayleigh quotient $v^t (\mathcal{P} + \mathcal{R}) v$, where $v \in \mathbb{S}^{n-1} \subset \mathbb{R}^n \backslash \{ 0 \}$.  The symmetry imposed on the curvature tensor implies that $\mathcal{P}+\mathcal{R} = \Sigma  + \Lambda$, where $\Sigma := \text{diag}(2R_{ k\overline{k} k \overline{k}})$ and $\Lambda$ is anti-symmetric.  The real bisectional curvature in any fixed frame is the Rayleigh quotient $v^t \mathcal{R} v$. The same argument shows that $v^t \mathcal{R} v = \Sigma + \Lambda'$, where $\Sigma$ is as before, and $\Lambda'$ is an anti-symmetric matrix. The Rayleigh quotient annihilates the anti-symmetric part, hence, $$\text{HSC}_{\omega} = v^t (\mathcal{P} + \mathcal{R}) v = v^t( \Sigma + \Lambda )v = v^t \Sigma v = v^t(\Sigma + \Lambda' ) v = \text{RBC}_{\omega}.$$
\end{proof}

\section{The Hermitian Chern--Lu Inequality}
For the sake of completeness, let us recall the proof of the Chern--Lu inequality given in \cite{WuYau1, YangZhengRBC}. Let $f : (M, \omega) \longrightarrow (N, \eta)$ be a holomorphic map between Hermitian manifolds. The Bochner formula gives us a formula for the complex Hessian of the square of the norm of a section $\sigma$ of a holomorphic vector bundle $\mathcal{E}$  \begin{eqnarray}\label{Bochner1}
\partial \overline{\partial} | \sigma |^2 &=& \langle \nabla^{1,0} \sigma, \nabla^{1,0} \sigma \rangle -  \langle R^{\mathcal{E}} \sigma, \sigma \rangle,
\end{eqnarray}

where $R^{\mathcal{E}}$ is the curvature of the Chern connection $\nabla$ on $\mathcal{E}$. We will apply the Bochner formula to $\sigma = \partial f$, where $\mathcal{E} = \mathcal{T}_M^{\ast} \otimes f^{\ast} \mathcal{T}_N$. From the splitting of the curvature of the Chern connection, we know that \begin{eqnarray}\label{ChernSplit2}
R^{\mathcal{T}_M^{\ast} \otimes f^{\ast} \mathcal{T}_N} &=& -R^{\mathcal{T}_M} \otimes \text{id} + \text{id} \otimes f^{\ast} R^{\mathcal{T}_N}.
\end{eqnarray}

The Chern--Lu inequality comes from computing $\Delta_{\omega} | \partial f |^2$. Since $\Delta_{\omega} = \text{tr}_{\omega} \sqrt{-1} \partial \overline{\partial}$, taking the trace of \eqref{Bochner1} with respect to $\omega$, and using \eqref{ChernSplit2}, we have (using the invariant notation\footnote{Namely, $\omega^{\sharp}$ denotes the metric $g^{-1}$ on $\mathcal{T}_M^{\ast}$. } employed in \cite{RubinsteinKE}): \begin{eqnarray*}
\Delta_{\omega} | \partial f |^2 & = & | \nabla^{1,0} \partial f |^2 + \left( \text{Ric}_{\omega}^{(2)} \right)^{\sharp} \otimes \eta (\partial f, \overline{\partial f}) \\
&& \hspace*{3cm} - \omega^{\sharp} \otimes \omega^{\sharp} \otimes R^{\mathcal{T}_N} (\partial f, \overline{\partial f}, \partial f, \overline{\partial f}).
\end{eqnarray*}

In coordinates, the above formula reads: \begin{eqnarray}\label{CALC}
g^{i \overline{j}} \partial_i \partial_{\overline{j}} (g^{k \overline{\ell}} h_{\gamma \overline{\delta}} f_k^{\gamma} \overline{f_{\ell}^{\delta}}) &=& g^{i \overline{j}} g^{k \overline{\ell}} h_{\gamma \overline{\delta}} f_{ik}^{\gamma} \overline{f_{j\ell}^{\delta}} + g^{i \overline{j}} R_{i \overline{j} k \overline{\ell} }^g g^{k \overline{q}} g^{p \overline{\ell}} h_{\alpha \overline{\beta}} f_p^{\alpha} \overline{f_q^{\beta}} \nonumber \\
&& \hspace*{4cm} - g^{i \overline{j}} g^{p \overline{q}} R_{ \alpha \overline{\beta}\gamma \overline{\delta}}^h f_i^{\alpha} \overline{f_{j}^{\beta}} f_p^{\gamma} \overline{f_{q}^{\delta}}.
\end{eqnarray}

\subsection*{Reminder 3.1}\label{SecondRicci}
For readers not so familiar with non-K\"ahler curvature tensors, let us recall that the trace of the curvature of the Chern connection yields three Chern Ricci curvatures: \begin{eqnarray*}
\text{Ric}^{(1)} \ : = \ g^{k \overline{\ell}} R_{i \overline{j} k \overline{\ell}}, \hspace*{1cm} \text{Ric}^{(2)} \ : = \ g^{i \overline{j}} R_{i \overline{j} k \overline{\ell}}, \hspace*{1cm} \text{Ric}^{(3)} \ : = \ g^{i \overline{\ell}} R_{i \overline{j} k \overline{\ell}}.
\end{eqnarray*}

Note that $\text{Ric}^{(4)} := g^{k \overline{j}} R_{i \overline{j} k \overline{\ell}}$ is merely the conjugate of $\text{Ric}^{(3)}$, inherited from the symmetry $\overline{R_{i \overline{j} k \overline{\ell}}} = R_{j \overline{i} \ell \overline{k}}$.\\

The following is due to \cite{YangZhengRBC}:

\subsection*{Lemma 3.2}
Let $f : (M, \omega) \longrightarrow (N, \eta)$ be a holomorphic map between Hermitian manifolds. Then away from the critical points of $f$, we have \begin{eqnarray*}
| \partial f |^2 \Delta_{\omega} \log | \partial f |^2 & \geq &  \text{Ric}_{\omega}^{(2)}  \otimes \eta (\partial f, \overline{\partial f}) - \omega^{\sharp} \otimes \omega^{\sharp} \otimes R^{\mathcal{T}_N} (\partial f, \overline{\partial f}, \partial f, \overline{\partial f}).
\end{eqnarray*}

In coordinates, this reads \begin{eqnarray*}
| \partial f |^2 \Delta_{\omega} \log |\partial f |^2 & \geq & \text{Ric}_{k \overline{\ell}}^{g} g^{k \overline{q}} g^{p \overline{\ell}} h_{\alpha \overline{\beta}} f_p^{\alpha} \overline{f_q^{\beta}} - g^{i \overline{j}} g^{p \overline{q}} R_{ \alpha \overline{\beta}\gamma \overline{\delta}}^h f_i^{\alpha} \overline{f_{j}^{\beta}} f_p^{\gamma} \overline{f_{q}^{\delta}}.
\end{eqnarray*}

\begin{proof}
For a smooth positive function $u : M \longrightarrow \mathbb{R}$, we have \begin{eqnarray*}
\Delta_{\omega} \log(u) &=& \frac{\Delta_{\omega} u}{u} - \frac{| \partial u |_{\omega}^2}{u^2}.
\end{eqnarray*}
Taking $u = | \partial f |^2$, we observe that $| \partial u |^2 = | \partial | \partial f |^2 |^2 = 4 | \partial | \partial f | |^2 | \partial f |^2 = | \nabla | \partial f | |^2 | \partial f |^2$. In particular, from \eqref{CALC} we have \begin{eqnarray*}
| \partial f |^2 \Delta_{\omega} \log | \partial f |^2 &=& \Delta_{\omega} | \partial f |^2 - | \nabla | \partial f | |^2 \\
&=& | \nabla  \partial f  |^2 + \text{Ric}_{k \overline{\ell}}^{g} g^{k \overline{q}} g^{p \overline{\ell}} h_{\alpha \overline{\beta}} f_p^{\alpha} \overline{f_q^{\beta}} - g^{i \overline{j}} g^{p \overline{q}} R_{ \alpha \overline{\beta}\gamma \overline{\delta}}^h f_i^{\alpha} \overline{f_{j}^{\beta}} f_p^{\gamma} \overline{f_{q}^{\delta}} - | \nabla | \partial f | |^2.
\end{eqnarray*}
The Kato inequality \cite{KatoReference} gives $| \nabla | \partial f | | \leq | \nabla \partial f |$, and hence \begin{eqnarray*}
| \partial f |^2 \Delta_{\omega} \log | \partial f |^2 & \geq &\text{Ric}_{k \overline{\ell}}^{g} g^{k \overline{q}} g^{p \overline{\ell}} h_{\alpha \overline{\beta}} f_p^{\alpha} \overline{f_q^{\beta}} - g^{i \overline{j}} g^{p \overline{q}} R_{ \alpha \overline{\beta}\gamma \overline{\delta}}^h f_i^{\alpha} \overline{f_{j}^{\beta}} f_p^{\gamma} \overline{f_{q}^{\delta}}.
\end{eqnarray*}
\end{proof}


We now complete the proof of \nameref{HermitianCL}. 
\begin{proof}
From \nameref{SecondRicci}, the curvature term $g^{i \overline{j}} R_{i \overline{j} k \overline{\ell}}^g$ is the second Ricci curvature $\text{Ric}_g^{(2)}$. From the assumption on the Ricci curvature of $\omega$, we have \begin{eqnarray}\label{RICLOWER1}
\text{Ric}_{k \overline{\ell}}^{g} g^{k \overline{q}} g^{p \overline{\ell}} h_{\alpha \overline{\beta}} f_p^{\alpha} \overline{f_q^{\beta}} & \geq & - C_1 g^{p \overline{q}} h_{\alpha \overline{\beta}} f_p^{\alpha} \overline{f_q^{\beta}} + C_2 g^{k \overline{q}} g^{p \overline{\ell}} h_{\alpha \overline{\beta}} h_{\gamma \overline{\delta}} f_p^{\alpha} \overline{f_q^{\beta}} f_k^{\gamma} \overline{f_{\ell}^{\delta}}.
\end{eqnarray}

The first term on the right-hand side of \eqref{RICLOWER1} is clearly $-C_1 \text{tr}_{\omega}(f^{\ast} \eta)$. For the second term, choose coordinates such that the metrics are Euclidean and $f_i^{\alpha} = \lambda_i \delta_i^{\alpha}$, where $\lambda_1 \geq \lambda_2 \geq \cdots \lambda_r \geq \lambda_{r+1} \geq 0$, where $r$ is the rank of $\partial f = (f_i^{\alpha})$. In these coordinates, we see that \begin{eqnarray*}
g^{k \overline{q}} g^{p \overline{\ell}} h_{\alpha \overline{\beta}} h_{\gamma \overline{\delta}} f_p^{\alpha} \overline{f_q^{\beta}} f_k^{\gamma} \overline{f_{\ell}^{\delta}} &=& \delta_k^q \delta_p^{\ell} \delta_{\alpha}^{\beta} \delta_{\gamma}^{\delta} \lambda_p \delta_p^{\alpha} \lambda_q \delta_q^{\beta} \lambda_k \delta_k^{\gamma} \lambda_{\ell} \delta_{\ell}^{\delta} \ = \ \lambda_{\alpha}^4.
\end{eqnarray*}

Since $\left( \sum_{\alpha} \lambda_{\alpha}^2 \right)^2 \leq r \sum_{\alpha} \lambda_{\alpha}^4$, we see that \begin{eqnarray*}
\text{Ric}_{k \overline{\ell}}^{g} g^{k \overline{q}} g^{p \overline{\ell}} h_{\alpha \overline{\beta}} f_p^{\alpha} \overline{f_q^{\beta}} & \geq & - C_1 | \partial f |^2 + \frac{C_2}{r} |\partial f |^4.
\end{eqnarray*}

These coordinates will also be useful for analyzing the remaining curvature term. Indeed, in these coordinates, we see that  \begin{eqnarray*}
g^{k \overline{\ell}} g^{i \overline{j}} R_{\gamma \overline{\delta} \alpha \overline{\beta}}^h f_i^{\alpha} \overline{f_j^{\beta}} f_k^{\gamma} \overline{f_{\ell}^{\delta}} &=& \delta_k^{\ell} \delta_i^j R_{\gamma \overline{\delta} \alpha \overline{\beta}}^h \lambda_i \delta_i^{\alpha} \lambda_j \delta_j^{\beta} \lambda_k \delta_k^{\gamma} \lambda_{\ell} \delta_{\ell}^{\delta} \ = \ R_{\gamma \overline{\gamma} \alpha \overline{\alpha}}^h \lambda_{\gamma}^2 \lambda_{\alpha}^2.
\end{eqnarray*}

Since $\text{RBC}_{\eta} \leq - \kappa \leq 0$, we have \begin{eqnarray*}
R_{\alpha \overline{\alpha} \gamma \overline{\gamma}}^h \lambda_{\alpha}^2 \lambda_{\gamma}^2 & \leq & - \kappa \sum_{\alpha} \lambda_{\alpha}^4 \ \leq \ - \frac{\kappa}{r} \left( \sum_{\alpha} \lambda_{\alpha}^2 \right)^2 \ = \ - \frac{\kappa}{r} |\partial f |^4.
\end{eqnarray*}
Combining these estimates yields \begin{eqnarray*}
\Delta_{\omega} \log | \partial f |^2 & \geq & - C_1 + \frac{1}{r} \left( C_2 + f^{\ast} \kappa \right) | \partial f |^2.
\end{eqnarray*}
If $M$ is compact,  and the functions $C_1, C_2, \kappa$ are constant, then the maximum principle yields the estimate $$|\partial f |^2 \ \leq \ \frac{C_1 r}{C_2 + \kappa}.$$ 
\end{proof}

\section{The Hermitian Aubin--Yau Inequality and Unification of the Schwarz Lemmas}
Let us now give a proof of the Hermitian Aubin--Yau second-order estimate \nameref{HermitianAubinYau}:

\begin{proof}
Let $f : (M,\omega) \to (N,\eta)$ be a holomorphic map, biholomorphic onto its image.  We saw before that \begin{eqnarray*}
\partial \overline{\partial} | \partial f |^2 &=& \langle \nabla^{1,0} \partial f, \nabla^{1,0} \partial f \rangle - \langle R^{\mathcal{T}_M^{\ast} \otimes f^{\ast} \mathcal{T}_N} \partial f, \partial f \rangle.
\end{eqnarray*}

Taking the trace with respect to $\eta$, we see that \begin{eqnarray*}
\Delta_{\eta} | \partial f |^2 &=& | \nabla^{1,0} \partial f |^2 - \omega^{\sharp} \otimes \text{Ric}_{\eta}^{(2)}(\partial f, \overline{\partial f}) \\
&& \hspace*{3cm} +  \left( R^{\mathcal{T}_M} \right)^{\sharp} \otimes   \eta \otimes \eta^{\sharp} (\partial f, \overline{\partial f}, \partial f^{-1}, \overline{\partial f^{-1}}),
\end{eqnarray*}

where $(R^{\mathcal{T}_M})^{\sharp}$ is of type $(2,2)$. In coordinates, this reads, \begin{eqnarray*}
h^{\gamma \overline{\delta}} \partial_{\gamma} \partial_{\overline{\delta}} \left( g^{i \overline{j}} h_{\alpha \overline{\beta}} f_i^{\alpha} \overline{f_j^{\beta}} \right) &=& h^{\gamma \overline{\delta}} g^{i \overline{j}} h_{\alpha \overline{\beta}} f_{ik}^{\alpha} \overline{f_{j\ell}^{\beta}} (f^{-1})_{\gamma}^k \overline{(f^{-1})_{\delta}^{\ell}} - h^{\gamma \overline{\delta}} g^{i \overline{j}} R^h_{\gamma \overline{\delta} \alpha \overline{\beta}} f_i^{\alpha} \overline{f_j^{\beta}} \\
&& \hspace*{3cm} + h^{\gamma \overline{\delta}} g^{i \overline{q}} g^{p \overline{j}} R^g_{k \overline{\ell} p \overline{q}} h_{\alpha \overline{\beta}} f_i^{\alpha} \overline{f_j^{\beta}} (f^{-1})_{\gamma}^k \overline{(f^{-1})_{\delta}^{\ell}}.
\end{eqnarray*}

Hence, from the upper bound on the second Ricci curvature, we have \begin{eqnarray*}
- \text{Ric}^h_{\alpha \overline{\beta}} g^{i \overline{j}} f_i^{\alpha}  \overline{f_j^{\beta}} & \geq & C_1 h_{\alpha \overline{\beta}} g^{i \overline{j}} f_i^{\alpha} \overline{f_j^{\beta}} - C_2 g_{k \overline{\ell}} g^{i \overline{j}} f_i^{\alpha} \overline{f_j^{\beta}} (f^{-1})_{\alpha}^k \overline{(f^{-1})_{\beta}^{\ell}} \\
&=& C_1 | \partial f |^2 - C_2 n
\end{eqnarray*}

Choose coordinates such that $g_{i \overline{j}} = \delta_{ij}$, $h_{\alpha \overline{\beta}} = \delta_{\alpha \beta}$, and $f_i^{\alpha} = \lambda_i \delta_i^{\alpha}$. Then \begin{eqnarray*}
h^{\gamma \overline{\delta}} g^{i \overline{q}} g^{p \overline{j}} R^g_{k \overline{\ell} p \overline{q}} h_{\alpha \overline{\beta}} f_i^{\alpha} \overline{f_j^{\beta}} (f^{-1})_{\gamma}^k \overline{(f^{-1})_{\delta}^{\ell}} &=& R_{k \overline{k} i \overline{i}}  \lambda_i^2 \lambda_k^{-2}.
\end{eqnarray*}

Assuming a lower bound on the Schwarz bisectional curvature, we have \begin{eqnarray*}
\sum_{i,k=1}^n R_{i \overline{i} k \overline{k}}  \lambda_i^2 \lambda_k^{-2} & \geq & - \kappa.
\end{eqnarray*}

Combining these estimates, we have \begin{eqnarray*}
\Delta_{\eta} | \partial f |^2 & \geq & C_1 | \partial f |^2 - n C_2 -  \kappa.
\end{eqnarray*}
\end{proof}

Applying this to the identity map, we get the following immediate corollary: 

\subsection*{Corollary 4.1} 
Let $M^n$ be a compact complex manifold that supports Hermitian metrics $\omega$, $\eta$. Assume that $\text{SBC}_{\omega} \geq - \kappa$ for some constant $\kappa \geq 0$, and suppose $\text{Ric}_{\eta}^{(2)} \leq - C_1\eta + C_2 \omega$ for some constants $C_1 > 0$ and $C_2 \in \mathbb{R}$. Then \begin{eqnarray*}
\text{tr}_{\omega}(\eta) & \leq & \frac{\kappa+ n C_2}{C_1}.
\end{eqnarray*}

\subsection*{Corollary 4.2}
Let $M$ be a compact complex manifold which supports Hermitian metrics $\omega$, $\eta$ with $\text{SBC}_{\omega} \geq - \kappa$ and $\text{Ric}_{\eta}^{(2)} \leq - C_1\eta  + C_2 \omega$, for constants $\kappa, C_1, C_2$, with $C_1>0$. If $\kappa \leq - C_2$, then the automorphism group $\text{Aut}(M)$ is trivial.\\

We can combine the Chern--Lu and Aubin--Yau inequalities proved in the previous section to obtain an $8$ $\mathbb{R}$--dimensional family of Schwarz lemmas:

\subsection*{Theorem 4.3}
Let $f : (M^n, \omega) \longrightarrow (N, \eta)$ be a holomorphic map between Hermitian manifolds with $\text{SBC}_{\omega} \geq - \kappa_1$ and $\text{RBC}_{\eta} \leq - \kappa_2$, for some constants $\kappa_1, \kappa_2 \geq 0$. Assume there is a Hermitian metric $\mu$ on $M$ such that, for constants $C_1, C_2, C_3, C_4 \in \mathbb{R}$, with $C_2 \geq 0$,  $C_3>0$, $\kappa_2 + C_2 >0$, we have \begin{eqnarray*}
- C_1 \mu + C_2 f^{\ast} \eta \ \leq \ \text{Ric}^{(2)}_{\mu} \ \leq \ - C_3 \mu + C_4 \omega.
\end{eqnarray*}

Then, if $M$ is compact, \begin{eqnarray*}
| \partial f|^2 & \leq & \frac{C_1 n r(\kappa_1 + C_4)}{C_3(\kappa_2 + C_2)},
\end{eqnarray*}

where $r$ is the rank of $\partial f$.\\

One particular corollary of the above theorem is the following Hermitian analog of the Chen--Cheng--Lu \cite{ChenChengLu} Schwarz lemma: 

\subsection*{Corollary 4.4}
Let $f : (M^n, \omega) \longrightarrow (N, \eta)$ be a holomorphic map of rank $r$ between Hermitian manifolds with $\text{SBC}_{\omega} \geq - \kappa_1$ and $\text{RBC}_{\eta} \leq - \kappa_2$, for some constants $\kappa_1, \kappa_2 \geq 0$. Assume there is a Hermitian metric $\mu$ on $M$ such that $$-C_1 \mu + C_2 f^{\ast} \eta \ \leq \  \text{Ric}_{\mu}^{(2)} \ \leq \ - C_3 \mu,$$ where $C_1 = \frac{\kappa_2 + C_2}{\kappa_2 nr} C_3$, and $C_2 \geq \kappa_2(nr-1)$. Then, if $M$ is compact, \begin{eqnarray*}
|\partial f|^2 & \leq & \frac{\kappa_1}{\kappa_2}.
\end{eqnarray*}

The theorem also yields a Schwarz lemma expressed exclusively in terms of (second) Ricci curvatures: 

\subsection*{Corollary 4.5}
Let $f : (M^n, \omega) \longrightarrow (N, \eta)$ be a holomorphic map of rank $r$ between Hermitian manifolds with $\text{SBC}_{\omega} \geq \kappa_1$ and $\text{RBC}_{\eta} \leq - \kappa_2$, for some constants $\kappa_1, \kappa_2 \geq 0$. Assume there is a Hermitian metric $\mu$ on $M$ such that $$-C_1 \mu + C_2 f^{\ast} \eta \ \leq \  \text{Ric}_{\mu}^{(2)} \ \leq \ - C_3 \mu + C_4 \omega,$$ where $C_2 \geq 0, C_3 >0$ and $nr(\kappa_1 + C_4) \leq \kappa_2 + C_2$. Then if $M$ is compact, $$| \partial f |^2 \ \leq \ \frac{C_1}{C_3}.$$ 

Of course, one can generate other Schwarz lemmas of Chen--Cheng--Lu-type by varying the constants appropriately. 

\subsection*{Non-compact source manifold}
In all the present formulations of the Schwarz lemma, we have given a formula for the Laplacian of the energy density, and then deduced an estimate on $| \partial f |^2$ under the assumption that the source manifold is compact. This is, of course, so that we can ensure the existence of a maximum point.  The compactness assumption can be relaxed to completeness of the source metric if the Riemannian Ricci curvature of the source metric is bounded below, and the (Chern) torsion $1$--form is uniformly bounded (see \cite[p. 13]{YangZhengRBC} for details).

\hfill

\scshape{Mathematical Sciences Institute, Australian National University, Acton, ACT 2601, Australia}

\scshape{BICMR, Peking University, Beijing, 100871, People's republic of china}

\textit{E-mail address}: \texttt{kyle.broder{@}anu.edu.au}

\end{document}